\input AHTOH-E.STY
\newif\ifanon
\ifanon\else\let\censored\relax\fi

\UDC{
512.772.7   
+ 512.624.3 
+ 512.552   
}


\title{%
Balanced factorisations
}
\author{%
\censored{
Anton A. Klyachko}$^\sharp$
\quad
and
\quad
\censored{Anton N. Vassilyev}$^\flat$
}
\address{
\ifanon
\else
\strut
$^\sharp$Faculty of Mechanics and Mathematics of Moscow State University
\\
Moscow 119991, Leninskie gory, MSU\\
klyachko@mech.math.msu.su
\\
$^\flat$Kazakhstan branch of
Moscow State University
\\
Astana 010010, ul. Kazhimukana 11, MSU
\\
antonvassilyev@mail.ru
\\
\fi
}

\grantsFirst{\ifanon{
\small
some grant kept secret for the
purpose of the double-blind review process}\else\RFBR15-01-05823\fi}
\grantsSecond{\ifanon
\small
another
grant
also kept secret%
\else Science committee of Ministry of Education and Science of
Republic of Kazakhstan, project no.  $\Gamma\Phi$4-0816\fi}

\abstract{%
{\it Any rational number can be factored into a product of several
rationals whose sum vanishes.} This simple but nontrivial fact was
suggested as a problem on a mathematical olympiad for high school
students.
We completely solve similar questions in all
finite fields and in some other rings, e.g., in the complex and real
matrix algebras. Also, we state several open questions.}

\s 0.
Introduction

\vskip-6mm


\disp{%
{\it ``Prove that
any rational number can be factored into a product of several
rationals whose sum vanishes."}}

\vskip-3mm

This problem
was
\ifanon\else 
invented by the second author and
\fi
suggested at the Kazakhstan republican mathematical olympiad for
high-school students in~2013 [Vas13]. A similar question about arbitrary
fields of characteristic not two was suggested at the Algebra olympiad
for university students at Moscow university in 2014 [Vas14]. Afterwards,
we learnt that the problem had been considered earlier [Iva13] (also in an
educational context).

The existence of such balanced factorisations is easy to prove in any
field of characteristic not two (see Theorem~0 below).  However, the
question on the possible numbers of factors in such factorisations is much
more difficult. This question is the main subject of our paper. For
example, any rational admits a balanced factoring into a product of five
factors, but some rationals do not admit balanced factorings into
products of three factors [Iva13]; the question about four factors is open
and seems to be difficult.%
\fn{%
When this paper was written, we learned that this question has a positive
answer [KMP16].
}
For instance, the author of
[Iva13] reproduced the following letter by M.~A.~Tsfasman
to him:
$$
\hbox{\framed{
$
\eqalign{\cr
3=(363/70)\cdot(20/77)\cdot(-49/110)\cdot(-5). \qbox{
\tencyrit Uf\dots
}&
\cr
\hbox{
\tencyrit Vash M.A.
}&
\cr\cr}
$
}}
$$
This is in Russian but no translation is needed --- the letter contains
the first discovered balanced decomposition of~3 into four
factors in the field of rationals (along with an interjection and
signature).  Such decompositions of
1 and~2
are less impressive:
$1=1\cdot1\cdot(-1)\cdot(-1)$
and
$2={1\over6}\cdot{9\over2}\cdot(-{2\over3})\cdot(-4)$.
[Iva13] contains computer-generated balanced
decompositions of first fifty positive integers into products of four
rational factors.

A similar problem for finite fields seems to be easier. Indeed, in each
\spacing{given} finite field, we can use a brute-force search and find all
element admitting balanced decompositions into any
\spacing{given} number of factors. This is what we actually did until
we realised that some more advanced algebra gives a
complete and computer-free solution to the problem in all finite fields.

One of our
main results (Theorem 2) describes all pairs $(q,k)$ such that
every element of the $q$-element field $\F_q$ admits a balanced
decomposition into a product of~$k$ factors.
The answer is nontrivial and rather
complicated. For instance, it turns out that, in all finite
fields except exactly one, each element admits a balanced factoring into a
product of at most three factors. The role of the unique exception is
played by the seven-element field~$\F_7$.  The main tool of our study of
finite fields is Hasse's estimate of the number of rational points of an
elliptic curve over a finite field.

In Section 2, we prove also that, in each field of characteristic
not 2, there is a ``universal" formula allowing us to obtain a balanced
factorisation of almost any element. For example, formula (1) gives a
balanced factorisation into five factors for any nonzero element (in any
field of characteristic not two). We prove that similar formulae exist
for~6,~7, and any larger numbers of factors but do not exist for three
factors. This fact is deduced from the Mason--Stothers theorem (the
$abc$-theorem for polynomials).

In Section 3, we show that the question about balanced factorisations in
finite-dimensional algebras is in essence reduced to a similar question
about fields. This allows us to solve the problem completely for some
natural algebras, e.g., the matrix algebras over $\C$ and $\R$.
The last section contains a list of questions remaining open.

The authors thank Yu. G. Prokhorov, M. A. Tsfasman, and also
(undergraduate students) Evgenia Kosheleva, Alisa~Pikulina, Nadira
Shoketaeva, and (a high-school student) Rauan Zhakypbek for useful
discussions. We are also grateful to anonymous referees for valuable
remarks that allowed us to improve the paper.

\s 1.
General results and remarks

The formal definition of our main
object looks as follows: suppose that an element $a$ of a ring is factored
into a product:  $a=a_1a_2\dots a_k$; we call this factorisation
\emph{balanced} if $\sum a_i=0$.

First note that over algebraically closed fields
the problem is trivial:
\disp
{\sl for any $k\ge2$ any element of any algebraically closed field
admits a balanced decomposition into a product of $k$ factors.
}
Indeed, for a given $a$ and $k$, we should find $a_1,\dots,a_k$ such that
$a=a_1\dots a_k$ and $a_1+\dots+a_k=0$. The solution is
straightforward: take, for example, $a_3=\dots=a_k=1$ and find
$a_1$ and $a_2=2-k-a_1$ from the quadratic equation $a_1(2-k-a_1)=a$.

The second observation is that the problem is easy (though non-trivial)
for any field of characteristic not two provided the number of factors is
at least five.

\Th 1.
For each $k\ge5$, in any field of characteristic not two, every element
decomposes into a product of $k$ factors whose sum vanishes. For each
$k<5$, there exists a field of characteristic not two where a similar
assertion is false.

\Proof
Let us prove the first assertion.
For zero element, we have nothing to prove; for nonzero element $a$, we
can solemnly write
$$
a={a\over2}\cdot{a\over2}\cdot(-a)\cdot{2\over a}\cdot\(-{2\over a}\).
\eqno{(1)}
$$
This gives a balanced decomposition into five factors.
A slight modification of (1) gives a balanced decomposition
of any element $b$ into six factors:
$$
b=c^2-ca=
{a\over2}\cdot{a\over2}\cdot(c-a)\cdot{2\over a}\cdot\(-{2\over a}\)
\cdot(-c),
\qbox{where $c$ is an element such that
$0\ne c^2\ne b$, and $a={c^2-b\over c}$}.
$$
(Such $c$ exists, except the case where
the field is $\F_3$ and $b=1$;
in this exceptional case we can take the factoring
$b=1=1^3(-1)^3$.)

A balanced decomposition into $k\ge7$ factors can be obtained
by multiplying one of the above decompositions
and a balanced decompositions of minus
one into two factors:
$
-1=(-1)\cdot1.
$
For example, we obtain the following balanced decomposition of any
element into a product of 100 factors:
$$
-b=-(c^2-ca)=
{a\over2}\cdot{a\over2}\cdot(c-a)\cdot{2\over a}\cdot\(-{2\over a}\)
\cdot(-c)
\cdot
(-1)^{47}\cdot1^{47}.
$$
(Note that $-b$ is an arbitrary element if $b$ is an arbitrary element).
This completes the proof of the first assertion.

The second assertion follows from Theorem 2 (see the next section): for
$k\le3$, the field~$\F_7$ is a required example; for $k=4$, we can take
$\F_3$. This completes the proof of Theorem 1.

\medskip

Note that, in this study, we consider no decompositions as ``trivial".
We allow factors to be 1, or $-1$, or anything and the problem remains
non-trivial. Actually, the role of bad decompositions is played
by so-called \emph{power} decompositions, i.e. decompositions with all
factors equal. This is not \emph{a priori} clear why such factorisations
are useless, but look at Theorem~4.

\s 2.
Fields

\Th 2.
Suppose that $k\ge2$ is an integer and $F$ is a finite field.
Then,
in $F$, any element can be decomposed into a product of
$k$ factors whose sum vanishes if and only if
\itemitem{either}
$|F|=2$ and $k$ is even,
\itemitem{or}
$|F|=4$ and $k\ne3$,
\itemitem{or}
$|F|$ is a power of two but neither two nor four (and $k$ is arbitrary),
\itemitem{or}
$|F|\in\{3,5\}$ and $k\notin\{2,4\}$,
\itemitem{or}
$|F|=7$ and $k\notin\{2,3\}$,
\itemitem{or}
$|F|$ is neither a power of two nor three nor five nor seven and
$k\ne2$.
\enditem
{\rm In other words, the situation in finite fields is the following:}

\vfil\break

{
\tabskip1em plus 2em
\halign to \hsize{
#\hfil&\strut\hfil#\hfil&\hfil#\hfil&\hfil#\hfil&\hfil#\hfil&\hfil#\hfil
\cr
          &$k=2$&$k=3$&$k=4$&$k=5,7,9,\dots$&$k=6,8,10,\dots$
\cr
$\F_2$    & yes & no  & yes &     no        &  yes
\cr
$\F_3$    & no  & yes & no  &     yes       &  yes
\cr
$\F_4$    & yes & no  & yes &     yes       &  yes
\cr
$\F_5$    & no  & yes & no  &     yes       &  yes
\cr
$\F_7$    & no  & no  & yes &     yes       &  yes
\cr
$\F_8,
\F_{16},
\F_{32},
\F_{64},
\dots$    & yes & yes & yes &     yes       &  yes
\cr
$\F_9,
\F_{11},
\F_{13},
\F_{17},
\dots$    & no  & yes & yes &     yes       &  yes
\cr
}
}

\smallskip

\Proof
Let us arrange the proof by the columns of this table.

{\bf\noindent Case $k=2$.}
In a finite field of characteristic two, any element is a square
(because the order of the multiplicative group of such a field is odd),
i.e. each element is a product of two equal factors whose sum
vanishes, because the characteristic is two. If the characteristic of a
finite field is not two, then not every element is a square and,
therefore, not every element decomposes into a product of two factors
whose sum vanishes
(because the equality $a=x\cdot(-x)$ implies that $-a$ is a square;
so, if each element has a balanced decomposition into a product of two
factors, then each element is a square).

{\smallskip\bf\noindent Case $k=3$.}
If the characteristic is three, then the order of the multiplicative group
$q-1=3^k-1$ is not divisible by three and, therefore, each element
is a cube and the decomposition $a=bbb$ is as required (because $b+b+b=0$
in a field of characteristic three).

To study the fields of other characteristics, we
need the well-known Hasse's estimate (also known as the Hasse--Weil bound).

\proclaim{Hasse's estimate} {\rm(see, e.g., [Sil86], Theorem V.1.1)}.
The number of points of an elliptic curve (i.e. a nonsingular and
irreducible over the closure of the field projective curve of genus
one) over a finite $q$-element field $\F_q$ is at least $q+1-2\sqrt q$.
\newline
In particular, this is true for nonsingular and irreducible (over
the closure of the field) cubic curves in the projective plane over $\F_q$.

Let us continue the proof assuming that the characteristic
is not three. We have to show that the system of equations
$$
\cases{
x+y+z=0
\cr
xyz=a
}
\eqno{(2)}
$$
over a finite field $\F_q$ has at least one solution
for any $a\in\F_q$. In other words, we have to show that the cubic
(affine) curve defined by the equation
$$
xy(x+y)=-a
$$
has at least one point over $\F_q$.
In homogeneous coordinates, the corresponding
projective curve has the equation
$$
XY(X+Y)=-aZ^3,
\eqno{(3)}
$$
and singular points of this curve are the solutions of the
system of equations consisting of
equation (3) and its partial derivatives with respect to
$X$, $Y$, and $Z$:
$$
\cases{
XY(X+Y)=-aZ^3
\cr
2XY+Y^2=0
\cr
2XY+X^2=0
\cr
-3aZ^2=0
\cr
}.
\eqno{(4)}
$$
We assume that $a\ne0$, because if $a=0$, then system (2) has
an obvious solution (zero).
Therefore, (and since the characteristic is not three)
the last equation of (4) implies $Z=0$. The difference of the second and
third equations shows that $X=\pm Y$; now, the second equation shows that
$X$ and $Y$ are zero (recall that $\Char\F_q\ne3$). Thus, system (4)
has no nonzero solutions, i.e. our projective curve has no singular points
over the closure of the field (if $\Char\F_q\ne3$).  This automatically
implies that our curve is irreducible (and, therefore, elliptic), because
a reducible cubic curve always has a singular point (over the closure of
the field): this is a point of intersection of components.

Thus, we can apply Hasse's estimate and conclude that projective
cubic (3) has more than three points over the field $\F_q$
if the characteristic of this field is
not
three and
$q+1>2\sqrt q+3$.
This inequality holds for $q\ge8$.  Thus, for
$q\ge8$, the
projective curve contains more than
three points and, hence, the corresponding affine curve contains at least
one point, because the intersection of an irreducible cubic with
the line at infinity
cannot contain more than three points%
\fn{%
In the case
under consideration,
the curve contains precisely
three points at
infinity:  $(1,0,0)$, $(0,1,0)$, and $(1,-1,0)$ (in homogeneous
coordinates).}%
, i.e.  system (2) has a solution as required.

It remains to investigate the fields $\F_2$, $\F_4$, $\F_5$, and $\F_7$.

In $\F_2$, the unity obviously has no balanced decompositions into a
product of three factors
(because factors cannot be zero, but $1+1+1\ne1$).

In $\F_4$, any nonzero balanced product $xyz$ of three factors
cannot contain equal factors (because $x+x=0$), therefore,
there is exactly one such product: this is the product of all nonzero
elements of the field and it equals one; hence, elements different
from one and zero do not admit balanced decompositions into products of
three factors.

In $\F_5$, system (2) has a solution: $x=y=b$, $z=-2b$, where $b$ is a
cubic root of $-{a\over2}$ (in $\F_5$, any element is a cube).

The seven-element field indeed is an exception:
$\pm3$ have no balanced decompositions into
products of three factors: if
$$
\cases{
x+y+z=0
\cr
xyz=\pm3
},
$$
then $x$, $y$, and $z$ must be pairwise different.
Indeed, if $y=x$, then $z=-2x$ and $\mp3=2x^3$, but
$\mp3$ is not twice a cube (cubes in~$\F_7$
are 0 and $\pm1$). Certainly, no two from
$x$, $y$, and $z$ can be opposite. So, only one
possibility remains (up to signs and permutations):
\newline
$
x=\pm1,
\quad
y=\pm2,
\qqbox{and}
z=\pm3.
$
But the product of such three numbers is $\pm1$, not $\pm3$.
(The element 3 has a shorter balanced
decomposition: $3=2\cdot(-2)$ but $-3$ has no such factorisations.)

{\smallskip\bf\noindent Case $k=4$.}
Let us try to obtain a balanced decomposition of an element $a\in F$
into a product of four factors, where one factor is 1.
The following argument (up to a point)
are similar to the proof in the case $k=3$.
We want to show that the system of equations
$$
\cases{
x+y+z+1=0
\cr
xyz=a
}
\eqno{(2')}
$$
over a finite field $\F_q$ has at least one solution
for any $a\in\F_q$. In other words, we want show that the cubic
(affine) curve defined by the equation
$$
xy(x+y+1)=-a
$$
has at least one point over $\F_q$. In homogeneous coordinates, the
corresponding projective curve has the equation
$$
XY(X+Y+Z)=-aZ^3
\eqno{(3')}
$$
and the singular points of this curve are the solutions
of the
system
consisting
of equation $(3')$ and its
partial derivatives with respect to $X$, $Y$, and $Z$:
$$
\cases{
XY(X+Y+Z)=-aZ^3
\cr
2XY+Y^2+YZ=0
\cr
2XY+X^2+XZ=0
\cr
XY=-3aZ^2
\cr
}.
\eqno{(4')}
$$
The difference of the second and third
equations is $Y^2-X^2+Z(Y-X)=0$. Thus,
either $X+Y+Z=0$ or $X=Y$.

If $X+Y+Z=0$, then the first equation of $(4')$ gives $Z=0$. Now, the
last equation of $(4')$ gives $XY=0$, and, therefore all
unknowns vanish (because we assume that $X+Y+Z=0$).

If  $X=Y$, then the second equation of $(4')$ shows that
$3X^2+XZ=0$. Here, if $X=0$, then and $Y=0$ and, therefore,
$Z=0$ (from the first equation of $(4')$). If $X\ne 0$,
we obtain $3X+Z=0$. Then, the last equation
of $(4')$ gives $27a=-1$. Thus, if $27a\ne-1$
and $|F|\ge8$, then we can apply
Hasse's estimate and conclude that $a\in F$ has a balanced decomposition
into a product of four factors (one of which is 1). If $27a=-1$,
we have a balanced decomposition of $a$:
$$
-{1\over27}=\(-{1\over3}\)\cdot\(-{1\over3}\)\cdot\(-{1\over3}\)\cdot1.
$$
(Actually, if $27a=-1$, the curve is singular but the singular point
itself is not a point at infinity and, hence, gives a balanced
factorisation of $a$.)

It remains to consider small fields $F$ with $|F|<8$.
In $\F_2$ and in $\F_4$ (as well as a in any finite field
of characteristic
two) any element is the fourth power of another element and this gives
a
balanced decomposition into a product of four (equal) factors.

In $\F_3$, the only nonzero balanced
product of four factors is $1\cdot1\cdot(-1)\cdot(-1)$ and
it equals 1; therefore, $-1$ does not admit such decompositions.

In $\F_5$, a product of four nonzero factors can be
one of the following:
$$
\(\pm1\)\(\pm1\)\(\pm1\)\(\pm1\),
\quad
\(\pm1\)\(\pm1\)\(\pm1\)\(\pm2\),
\quad
\(\pm1\)\(\pm1\)\(\pm2\)\(\pm2\),
\quad
\(\pm1\)\(\pm2\)\(\pm2\)\(\pm2\),
\quad
\(\pm2\)\(\pm2\)\(\pm2\)\(\pm2\).
$$
The first, third, and fifth products equal $\pm1$, because all squares
equal $\pm1$. In the second and fourth product, there are
only two arrangements of signs making the sum of factors
zero:
$$
\eqalign{
&1\cdot1\cdot1\cdot2,
\quad
(-1)\cdot(-1)\cdot(-1)\cdot(-2),
\cr
&(-1)\cdot2\cdot2\cdot2,
\quad
1\cdot(-2)\cdot(-2)\cdot(-2).
}
$$
All these products equal two; therefore, $-2\in\F_5$ admits
no balanced decomposition into a product of four factors.

In $\F_7$, we find explicit balanced decompositions:
$$
0=0^4,
\quad
1=1^2\cdot(-1)^2,
\quad
-1=1\cdot1\cdot2\cdot3,
\quad
2=2^2\cdot(-2)^2,
\quad
-2=1\cdot(-2)^2\cdot3,
\quad
3=(-1)\cdot2\cdot3^2,
\quad
-3=(-1)^3\cdot3.
$$

{\smallskip\bf\noindent Case of even $k>4$.}
If each element $a$ has a balanced decomposition into a product of
$k$ factors: $a=a_1\dots a_k$, then each element has a
balanced decomposition into a product of $(k+2)$ factors:
$-a=a_1\dots a_k\cdot1\cdot(-1)$
(because $-a$ is an arbitrary element if $a$ is an arbitrary element).
Therefore, it suffices to prove
the assertion for $k=6$. Moreover, for all finite fields, except $\F_3$
and $\F_5$, the assertion is true, because we have constructed
a balanced decomposition of each element into a product of four
factors.

In $\F_3$, we have $0=0^6$,\quad $1=1^6$,\quad $-1=1^3\cdot(-1)^3$.
In $\F_5$, the balanced product
$x\cdot x\cdot(-2x)\cdot1\cdot1\cdot(-2)$ equals $-x^3$ which is
any element, because all elements are cubes.

{\smallskip\bf\noindent Case of odd $k>4$.}
The same induction as in the case of even large $k$ makes it possible
to reduce the problem to the case $k=5$. Moreover, for all finite fields,
except $\F_2$, $\F_4$, and $\F_7$, the assertion is true,
because we have already constructed a balanced decomposition of each
element into a product of three factors.

In $\F_7$, the desired decomposition exists by Theorem 0.
In $\F_4$, we can write $a=b^2xyz$, where $b$ is a square root of~$a$
and $x,y,z$ are all nonzero elements of the field (their product is
one and their sum is zero).  In $\F_2$,
there are no
balanced decompositions of 1
into products of odd number of factors
(because factor cannot be zero and the sum of an odd number of unities is
not zero).  This completes the proof.

\bigskip

Formula (1) can be considered as a ``universal formula"
making it possible to factorise balancedly almost any element of any field
of characteristic not two into a product of five factors
(where \emph{almost any} means any, except a finite number
of elements). Theorem 0 shows that such a universal
formula exists for each $k\ge5$. The proof of Theorem~0
gives explicit formulae:
$$
t=
\underbrace{
{t\over2}\cdot{t\over2}\cdot(-t)\cdot{2\over t}\cdot\(-{2\over t}\)
}_{5\ factors}
=%
\underbrace{
{1-t\over2}\cdot{1-t\over2}\cdot t\cdot{2\over1-t}\cdot{2\over t-1}
\cdot(-1)
}_{6\ factors}
=%
\underbrace{
\(-{t\over2}\)\cdot\(-{t\over2}\)\cdot t\cdot\(-{2\over t}\)\cdot{2\over t}
\cdot(-1)\cdot1
}_{7\ factors}
=\dots
$$
The following theorem shows that no ``universal formula"
for balanced decompositions into three factors exists
(a universal balanced decomposition into two factors do not exist either
for an obvious reason; the question about four factor remains open, see the
last section).

\Th 3.
For any field $F$, the element $t$ of the field of rational fractions
$F(t)$ does not admit a balanced decomposition into a product of three
factors.


\Proof
Assuming the contrary (and finding a common denominator), we obtain
the identity
$$
t^s={x(t)\over v(t)}\cdot{y(t)\over v(t)}\cdot{z(t)\over v(t)},
\qbox{where $x,y,z\in F[t]$ and $x+y+z=0$.}
$$
We have to show that $s\ne1$, but we prefer to prove a stronger fact:
\disp{\sl the above equalities imply that
$s$ is a multiple of 3.}
The polynomials $x$, $y$, and~$z$ can be a assumed to be
coprime, because the equality ${xyz=t^sv^3}$ shows that
an irreducible common divisor of $x$, $y$, and $z$ must
either divide $v$ or be $t$;
in both cases,
the equation can be cancelled.
In addition, we may assume that $v(0)\ne0$ (increasing
$s$ if needed).

Let us recall the well-known Mason--Stothers theorem
([May84], [Sto81]) that can be found in many books
(see, e.g., [Lang02]).  We prefer to use the version due
to Snyder, which works in any characteristic.

\proclaim{Mason--Stothers theorem {\rm(in the form of Snyder [Sny00])}}{}.
If three polynomials
${x,y,z\in F[t]}$ over a field $F$ are
coprime and $x+y+z=0$,
then either the degrees of all these polynomials are strictly less
than the number of
different roots of the product~$xyz$ in the algebraic closure of $F$
or all three derivatives $x'$, $y'$, and $z'$ vanish (as polynomials).

In the case under consideration,
$xyz=t^sv^3$ and the number of different roots of this
polynomial is at most $\deg v+1$; therefore, the Mason--Stothers theorem
says that either the degree of each of $x,y,z$ is at most the
degree of~$v$ or $x'=y'=z'=0$.

In the first case, $\deg(xyz)\le 3\deg v$
and, hence, $s=0$ (because $xyz=t^sv^3$) as required.
In the second case, the derivative of the product vanishes:
$0=(xyz)'=(t^sv^3)'=st^{s-1}v^3+3t^sv^2v'=v^2t^{s-1}(sv+3tv')$;
cancelling $v^2t^{s-1}$, we obtain $sv=-3tv'$. This means that $s$
is divisible by $\Char F$, since $v(0)\ne0$.
Therefore, either $\Char F=3$ and $s$ is a multiple of three as required,
or $v'=0$.

If $v'=0$, let us recall that
an equality $f'=0$
means that the polynomial $f$ has the form
$$
f(x)=f_1(x^p),
\qbox{where $p$ is the characteristic of the field, and
$f_1$ is a polynomial.}
$$
Therefore, substituting
$$
x(t)=x_1(t^p), \quad y(t)=y_1(t^p), \quad z(t)=z_1(t^p), \quad
v(t)=v_1(t^p),
$$
to the initial identity,
we obtain that $s$ is divisible by $p$ and, putting $t^p=\tau$,
we arrive to a similar equality for polynomials of lower degree:
$$
\tau^{s/p}={x_1(\tau)\over v_1(\tau)}\cdot{y_1(\tau)\over v_1(\tau)}
\cdot
{z_1(\tau)\over
v_1(\tau)},
\qbox{where $x_1,y_1,z_1\in F[\tau]$ and $x_1+y_1+z_1=0$.}
$$
An obvious induction completes the proof.


\s 3.
Algebras

\Lemma 1.
Suppose that the value of
a one-variable polynomial over an associative commutative ring with
unity at some point $d$ is nilpotent
and the value of the derivative at this point is invertible.
Then the polynomial has a root in this ring.
Moreover, for some root $b$, the difference $d-b$ is divisible by $f(d)$.

\Proof
An obvious change of variables reduces the situation to case, where $d=0$.
Suppose that the polynomial over a ring $R$ has the form
$f(x)=a_0+a_1x+\dots+a_nx^n$, where $a_1$ is invertible
and $a_0^s=0$. We argue by induction on $s$ and have to prove that {\sl
$f$ has a root,divisible by $a_0$.}

In the quotient ring
$\=R=R/(a_0^{s-1}R)$, the image $\=f$ of
$f$ has a root $\=c\={a_0}$ by the induction hypothesis.
Take some preimage $c\in R$ of the element $\=c\in \=R$
and let us try to find a root of $f$ in the form $b=ca_0+ta_0^{s-1}$,
where $t$ is an (unknown) element of $R$. Since $a_0^s=0$, we have
$$
f(b)= a_0+a_1(ca_0+ta_0^{s-1})+\dots+a_n(ca_0+ta_0^{s-1})^n=
f(ca_0)+a_1ta_0^{s-1}.
\eqno{(5)}
$$
Now, $ca_0$ is a root of $f$ modulo the ideal
$a_0^{s-1}R$ and, hence, $f(ca_0)\in a_0^{s-1}R$, i.e.
$f(ca_0)=ra_0^{s-1}$ for some $r\in R$. It remains to note that,
in (5),
$f(b)$
vanishes if we take $t=-r/a_1$. This completes the
proof.

\smallskip

The following theorem reduces the question on balanced factorisations in
finite-dimensional algebras to a similar question in fields if
we take into account
only \emph{non-power} factorisations, i.e. factorisations
having at least two
non-equal factors.

\Th 4.
Let $F$ be a field and let $n$ be an integer larger than two. If,
in all finite extensions of~$F$, each element has a non-power balanced
decomposition into a product of $n$ elements, then the same is true for
each element of each finite-dimensional associative unital algebra
over~$F$.

\Proof
Clearly, it suffices to prove the assertion for
finite-dimensional
one-generator unital algebras (because any element of any algebra
lies in a one-generated subalgebra). Thus, we assume that
an algebra $A$ over $F$ has the form $A=F[x]/(f)$, where $f\in F[x]$.
Such algebra $A$ decomposes
into a direct sum
$$
A\iso\bigoplus_{i=1}^m F_i[x]/(x^{k_i}),
\qqbox{where fields $F_i$ are finite extensions of $F$}
$$
($F_i\iso F[x]/(p_i)$ if $f=\prod p_i^{k_i}$ is the decomposition
of $f$ into a product of irreducible (over $F$) factors).
It suffices to obtain a balanced decomposition for
each direct term. Therefore, we assume, that $A=G[x]/(x^k)$,
where the field $G$ is a finite extension of $F$. Such algebra $A$ is
local, i.e. it has a unique maximal ideal $I$ (generated by $x$),
$A/I\iso G$ and all elements not lying in $I$ are invertible.

We want to decompose any element $a\in A$ into a product of $n$
elements with zero sum.

{\noindent\bf Case I. $a\notin I$.}
In this case, we find a
non-power
balanced decomposition of $a$ modulo
ideal~$I$, i.e. in the field $G$. Thus, we obtain elements
$a_1,\dots,a_n\in A$ such that
$$
a-a_1a_2\dots a_n\in I,
\quad
a_1+\dots+a_n\in I
\qqbox{and (without loss of generality)}
a_1-a_n \notin I.
$$
This means that, for the quadratic polynomial
$$
g(t)=a+ta_2a_3\dots a_{n-1}(t+a_2+a_3+\dots+a_{n-1}),
\qqbox{we have}
g(a_1)\in I.
\eqno{(6)}
$$
For the derivative of $g$, we obtain
$$
g'(a_1)=a_2a_3\dots a_{n-1}(a_1+a_2+a_3+\dots+a_{n-1})+
a_1a_2a_3\dots a_{n-1}\in a_2a_3\dots a_{n-1}(a_1-a_n)+I.
$$
The ideal $I$ consists of nilpotent elements and all
elements of $A\setminus I$ are invertible. Therefore, the conditions of
Lemma~1 are satisfied, because $a_1\ne a_n \pmod I$. Applying Lemma 1, we
find a root $\~t\in A$ of~$g$ and obtain a decomposition:
$$
a=\~ta_2a_3\dots a_{n-1}(-\~t-a_2-a_3-\dots-a_{n-1})
\qbox{with zero sum of factors.}
\eqno{(7)}
$$
This decomposition is non-power, because $\~t\equiv a_1 \pmod I$
by Lemma 1 and $a_1\not\equiv a_n \pmod I$ by the assumption.

\medskip

{\noindent\bf Case II. $a\in I$.}
Let us choose an invertible (i.e. not lying in $I$)
elements $a_2,\dots,a_{n-1}\in A$
in such a way that their sum is also
invertible. This is possible if the field $G=A/I$ has more than
two elements. If the field $G$ is two-element, then
the unit element has in $G$ no non-power
decomposition that contradicts the condition.

For the polynomial $g(t)$ (see formula (6))
we obtain that $g(0)=a$ is a nilpotent element and
$$
{g'(0)= a_2a_3\dots a_{n-1}(a_2+a_3+\dots+a_{n-1})}
\hbox{ is an invertible element.}
$$
Therefore, by Lemma 1, $g$ has a root $\~t\in A$ as
required (see~(7)). Decomposition (7) cannot be
power, because $a_2$ is invertible but $a$ is not.
This completes the proof.

\Corollary 1.
Each element of a finite-dimensional associative unital algebra (over a
field) decomposes into a product~of
\itemitem{\rm a)}
three elements whose sum vanishes if the field is algebraically closed;
\itemitem{\rm b)}
five elements whose sum vanishes if the characteristic of the
field is not two.

\Proof
The first assertion follows immediately from Theorem 4, because, in
an algebraically closed field, each element has a non-power
balanced decomposition into a product of three factors (to
obtain a non-power balanced decomposition $a=a_1a_2a_3$ of
a given element $a$, we can choose any element
$a_1$ such that $a_1^3\ne a$ and then $a_2$ and $a_3$
can be found from a quadratic equation).

To prove the second assertion, it suffices to apply
Theorems 2 and 0
and note that formula~(1) always gives a non-power
decomposition.

\Corollary 2.
For any $k\ge3$,
any complex or real matrix can be decomposed into a product of
$k$ matrices (over the same field) whose sum vanishes.

\Proof
The assertion follows immediately from Theorem 4, because each
real or complex number $a$ admits a
nonpower balanced decomposition
$a=x\cdot(x+1)\cdot1^{k-3}\cdot(2-k-2x)$,
as this equality is a cubic equation
with respect to $x$.


Now, we give examples showing that no conditions of
Theorem 4 and its corollaries can be omitted.

\Example 1.
Each element of the tree-element field $\F_3$ has a
balanced decomposition into a product of three factors:
$
0=0\cdot0\cdot0,
\quad
1=1\cdot1\cdot1,
\quad
2=2\cdot2\cdot2.
$
However,
in the two-dimensional algebra $A=\F_3[x]/(x^2)$ over this field,
the element
$1+x$ does not admit balanced decompositions into a product of
three factors, because the decomposition
$1=1\cdot1\cdot1$ is the unique balanced decomposition
of 1 in $\F_3$; therefore, the balanced decomposition of $1+x\in A$
must have the form $1+x=(1+kx)(1+lx)(1+mx)$ (where $k,l,m\in\F_3$), whence
we obtain $k+l+m=1$ and the decomposition is not balanced. This example
shows that Theorem 4 become false if we omit the words \emph{non-power}.

\Example 2.
In the algebra of polynomials $F[x]$ over any field, the element $x$
has no balanced decompositions. This example shows that
finite-dimensionality condition cannot be omitted
in Theorem 4 and Corollary 1.

In algebras with zero multiplication, no nonzero element
has balanced decompositions. This shows that
the condition that the algebra has a unit also
cannot be omitted in Theorem 4
and Corollary 1.

The condition $n>2$ can be omitted in
Theorem 4, because this condition follows from
other conditions: in any field, any balanced decomposition of zero into a
product of two factors must be power.
On the other hand, in any nonzero ring, zero has non-power
decompositions into products of three and any larger numbers of
factors, e.g., $0=0^{\the\year}\cdot b\cdot(-b)$, where $b$ is a nonzero
element.  However, there is the following simple example.

\Example 3.
In the field of complex numbers, any \emph{nonzero} element has a
non-power balanced decomposition into a product of two factors,
but the nilpotent Jordan block obviously has no balanced
decomposition into a product of two factors (for any field),
because such a decomposition of $J$ would mean that $-J$ is a square, but
it is not.

Example 3 also shows that, in Corollary 2, we cannot omit the
condition $k>2$ and, in Corollary 1(a), it is impossible to replace
three with two. The following example shows that, in Corollary 1(b), we
cannot replace five with a lower number.

\Example 4.
As mentioned above (see Example 1), in the two-dimensional algebra
$A=\F_3[x]/(x^2)$, the element $1+x$ does not admit a balanced
decomposition into three factors. In the same algebra (as well as in the
field $\F_3$), minus one admits no balanced
decomposition into a product of four factors and 1 has no balanced
decompositions into products of two factors.

\Example 5.
In the field $\F_2$, the identity element does not admit balanced
decompositions into products of five factors. This simple example
shows that the condition on characteristic cannot be omitted in
Corollary 1(b) (and in Theorem~0).



\s 4.
Open questions

\Question 1 {\rm (A. V. Ivanishchuk [Iva13])}.
Can
any rational number be decomposed into a product
of four rational numbers whose sum vanishes?{\rm*$^)$}

\Question 2.
Can any element of any field be decomposed into a
product of at most four factors whose sum vanishes?

\Question 3.
Does there exist a universal formula for balanced decomposition
into four factor? More precisely, does the element $t$ of the field of
rational fractions~$\C(t)$ (or even $\Q(t)$)
admit a balanced decomposition into a product of four factors?%
\fn{%
When this paper was written, we learned that the answers to
Questions 1 and 3 are positive [KMP16]. }

\Question 4.
What does occur in characteristic 2\thinspace?
Does there exist universal
formulae? Does any element of any field admit a balanced
factorisation?


\REFERENCES

\[Vas13]
Vassilyev A. N.
Kazakhstan republican olympiad in mathematics. 2013.
Final stage. Ninth grade. Problem~4.
(in Russian)
{\tt http://matol.kz/olympiads/151}

\[Vas14]
Vassilyev A. N.
Ninth algebra olympiad for students in MSU. 2014. Problem 3.
(in Russian).
\newline
{\tt http://halgebra.math.msu.su/Olympiad/}

\[Iva13]
Ivanishchuk A. V.
The experience of learning and research activity of students
in Lyceum 1511 (MEPhI) (in~Russian)
//
published
in the book
{\it Sgibnev A. I.
Research problems for beginners.
Moscow: MCCME,~2013.}
(Freely available at {\tt http://www.mccme.ru/free-books/})

\[KMP16]
Klyachko A. A., Mazhuga A. M., Ponfilenko A. N.
Balanced factorisations in some algebras.
arXiv:1607.01957

\[Lang02]
Lang S.
Algebra.
New York, Berlin, Heidelberg: Springer-Verlag, 2002.

\[Mas84]
Mason R. C.,
Diophantine Equations over Function Fields,
London Mathematical Society Lecture Note Series~96,
Cambridge, England: Cambridge University Press, 1984.

\[Sil86]
Silverman J. H.
The Arithmetic of Elliptic Curves.
New York: Springer-Verlag, 1986.

\[Sny00]
Snyder N.
An alternate proof of Mason's theorem.
Elem. Math., 2000, 55:3, 93--94.

\[Sto81]
Stothers W. W.,
Polynomial identities and hauptmoduln,
Quarterly J. Math., 1981, 32:3, 349--370.

\end